\documentclass{amsart}

\usepackage{amsthm,amsmath,amssymb}
\usepackage{latexsym,hyperref}
\usepackage{tensor,mathdots}

\usepackage{tensor}
\usepackage{upgreek}

\theoremstyle{plain}
\newtheorem{theorem}{Theorem}[section]

\theoremstyle{definition}
\newtheorem{definition}[theorem]{Definition}

\theoremstyle{remark}
\newtheorem{remark} [theorem]{Remark}

\newtheorem{example}[theorem]{Example}
\author{Elvis Barakovic}
\address{
Department of Mathematics \newline \indent
Faculty of Science and Mathematics \newline \indent
University of Tuzla  \newline \indent
Univerzitetska 4 \newline \indent
75000 Tuzla \newline \indent
Bosnia and Herzegovina}
\email{elvis.barakovic@untz.ba}

\author{Vedad Pasic}
\address{
Department of Mathematics \newline \indent
Faculty of Science and Mathematics \newline \indent
University of Tuzla  \newline \indent
Univerzitetska 4 \newline \indent
75000 Tuzla \newline \indent
Bosnia and Herzegovina}
\email{vedad.pasic@untz.ba}
\urladdr{http://pmf.untz.ba/vedad/}

\title[Analyzing the spectral (a)symmetry of the massless Dirac operator]{Analyzing the spectral (a)symmetry of the massless Dirac operator on the 3-torus}


\begin{document}
\begin{abstract}
We analyze the spectrum of the massless Dirac operator on the 3-torus $\mathbb{T}^3$. It is known that it is possible to calculate this spectrum explicitly, that it is symmetric about zero and that each eigenvalue has even multiplicity. However, for a general oriented closed Riemannian 3-manifold $(M,g)$ there is no reason for the spectrum of the massless Dirac operator to be symmetric. Using perturbation theory, we derive the asymptotic formulae for its eigenvalues and prove that by the perturbation of the Euclidean metric on the 3-torus, it is possible to obtain spectral asymmetry of the massless Dirac operator in the axisymmetric case.
\end{abstract}
\maketitle

\section{Introduction}
We work on a $3$-dimensional connected compact oriented manifold $M$ equipped
with a Riemannian metric $g_{\alpha\beta}$. Our aim is to analyze the spectrum of the massless Dirac operator on $M$
which describes a single massless neutrino living in a 3-dimensional compact universe.
The relationship between the spectrum of an operator and the manifold geometry is of particular interest in recent research.
The eigenvalues of the massless Dirac operator represent the energy levels of the massless neutrino and we are in particular interested in studying the spectrum of that operator, i.e. the set of all eigenvalues of the operator. An advanced review of the theory of the Dirac operator in general can be found in e.g. \cite{lawson1989spin}.
It is however very difficult to determine the spectrum of the massless Dirac operator on an arbitrary manifold $M$ and
there are only two known examples where the spectrum of the massless Dirac operator can be  calculated explicitly: the unit 3-torus $\mathbb{T}^3$ equipped with Euclidean metric (see \cite{friedrich1984abhangigkeit}) and the unit 3-sphere $\mathbb{S}^3$ equipped with metric induced by the natural embedding of $\mathbb{S}^3$ in $\mathbb{R}^4$ (see \cite{bar1996dirac,trautman1995dirac}).

In this paper we choose to work on the unit 3-torus $\mathbb{T}^3$ parameterized by cyclic coordinates $x^\alpha$, $\alpha=1,2,3$ of period $2\pi$. Under the assumption that the metric is Euclidean, the massless Dirac operator corresponding to the standard spin structure reads
\begin{equation}\label{dirac operator with euclidean metric}
W=-i\left(\begin{array}{cc}
            \frac{\partial}{\partial x^3} & \frac{\partial}{\partial x^1}-i\frac{\partial}{\partial x^2} \\
            \frac{\partial}{\partial x^1}+i\frac{\partial}{\partial x^2} & -\frac{\partial}{\partial x^3}
          \end{array}
\right),
\end{equation} see e.g. \cite{downes2013spectral}.
Note that we denote the massless Dirac operator by $W$ as it is often also referred to as the \emph{Weyl} operator, especially in theoretical physics, whereas we previously examined in some detail solutions of Einstein-Weyl theory, which incorporates the massless Dirac equation and operator (see \cite{barakovic2017physical, pasic2009new, pasic2010new, pasic2014pp, pasic2015torsion, pasic2017axial, pasic2014new, pasic2005pp}), where we suggested spacetimes with torsion as metric-affine models for the massless neutrino, which in turn motivated our current work in spectral analysis.

The eigenvalues and the eigenfunctions of the operator (\ref{dirac operator with euclidean metric}) can be calculated explicitly. The spectrum of the massless Dirac operator on the unit torus $\mathbb{T}^3$ equipped with the Euclidean metric is as follows: zero is an eigenvalue of multiplicity two and for each $m\in \mathbb{Z}^3\setminus \{0\}$,  the eigenvalues are $\pm \|m\|$, which in general have even multiplicity. Even multiplicity of eigenvalues is a consequence of the fact that the massless Dirac operator (\ref{dirac operator with euclidean metric}) and the operator of charge conjugation
\begin{equation}\label{eigenvector Cv2}
C\left(
   \begin{array}{r}
     v_1 \\
     v_2 \\
   \end{array}
 \right)=
 \left(
   \begin{array}{r}
     -\overline{v_2} \\
    \overline{ v_1} \\
   \end{array}
 \right).
\end{equation}
commute, see \cite{downes2013spectral}.

It is very important to note that the spectrum of the operator (\ref{dirac operator with euclidean metric}) is symmetric about zero, as was demonstrated in  \cite{bar1996dirac,chervova2013spectral,trautman1995dirac}.
However, as emphasized in \cite{atiyah1973spectral,atiyah1975spectral1,atiyah1975spectral2,atiyah1976spectral} for a general oriented Riemannian 3-manifold $(M,g)$ there is no justifiable reason for the spectrum of the massless Dirac operator to be symmetric, as it would imply that in these two cases there is no difference between the properties of the massless neutrino and the massless antineutrino. Therefore the primary objective of our study is to break the spectral symmetry of the massless Dirac operator.

\section*{Acknowledgements}
We are very grateful to Dmitri Vassiliev, without whose help, guidance and mentorship, this paper and our involvement in the subject area would not be possible.

Both authors were supported by the Ministry of Education and Science of the Federation of Bosnia and Herzegovina, within the research project ``Teorijska analiza spektra Diracovog operatora bez mase'' (number 05-39-2472-1/17).

The authors acknowledge the COST Action CA15117 (CANTATA), the COST action CA16104 (GWVerse) and the COST action CA18108.

\section{Perturbing the massless Dirac operator}\label{SectionPerturbation}
As we work on the 3--torus $\mathbb{T}^3$, which has trivial topology, spectral asymmetry will be obtained by perturbing the Euclidean metric itself, i.e. we consider a metric $g_{\alpha\beta}(x;\upepsilon)$, whose components are smooth functions of coordinates $x^\alpha, \alpha=1,2,3$ and small real parameter $\upepsilon$, which satisfies $g_{\alpha\beta}(x;0) = \delta_{\alpha\beta}$.
This is the approach that was similarly successfully  used in \cite{downes2013spectral} for the eigenvalue with smallest modulus $\lambda=0$, where spectral asymmetry was achieved and two particular families of Riemannian metrics were presented for which the eigenvalue with smallest modulus can be evaluated explicitly. The main goal of this paper is to derive the asymptotic formulae for all the other eigenvalues and to see under which conditions we can create spectral asymmetry. Note that the behavior of eigenvalues of the massless Dirac operator under perturbations of the metric was studied in \cite{bourguignon1992spineurs}, on a much more abstract level.

Our first goal is to write down explicitly the massless Dirac operator. To do this, we use the concepts of \emph{frame} and \emph{coframe}, whose differential geometrical definition and properties were given in \cite{chervova2013spectral}.
According to \cite{kirby1989topology,stiefel1935richtungsfelder}, a $3$-dimensional oriented manifold is parallelizable and consequently, there exist smooth real vector fields $e_j(x)$, $j=1,2,3$
that are linearly independent in every point $x$ of the manifold $M$, which we call a \emph{frame}. We can assume that the vector fields $e_j(x)$ are orthonormal, and if  not, the orthonormality can always be achieved using the Gram-Schmidt process. The coordinate components of the vector  $e_j(x)$ are $\tensor{e}{_j^\alpha}(x)$, $\alpha=1,2,3$, where the so-called anholonomic or \emph{frame index}, denoted by the Latin letter $j$, enumerates the vector field and the holonomic or \emph{tensor index}, denoted by Greek letter $\alpha$, enumerates their components. The \emph{coframe} is defined as the triple of covector fields $e^k(x)$, $k = 1, 2, 3$ and the coordinate components of the vector  $e^k(x)$ are $\tensor{e}{^k_\alpha}(x)$, $\alpha=1,2,3$, where
$\displaystyle e^k{}_\beta:=\delta^{kj}g_{\beta\gamma}e_j{}^\gamma.
$
The frame is completely determined by the coframe, and vice versa, by the relation
$
e_j{}^\alpha e^k{}_\alpha=\delta_j{}^k.
$

We consider the \emph{perturbed} coframe to be a smooth real-valued matrix function $e^j{}_\alpha(x;\upepsilon)$, $j,\alpha = 1,2,3$ satisfying the conditions
\[
g_{\alpha\beta}(x;\upepsilon) = \delta_{jk} e^j{}_{\alpha}(x;\upepsilon)e^k{}_\beta(x;\upepsilon), \quad
e^j{}_\alpha (x;0) = \delta^j{}_\alpha,
\] as was explained in \cite{downes2013spectral}. The perturbed frame is considered to be the smooth real-valued matrix-function $e_j{}^\alpha(x;\upepsilon), j,\alpha=1,2,3$ defined by
$$
e_j{}^\alpha(x;\upepsilon) e^k{}_\alpha(x;\upepsilon) = \delta_j{}^k.
$$

\begin{remark}\label{RemarkAxisymmetric}
In this paper we choose the frame and the coframe so that
they depend on the coordinate $x^1$ only, so that the original eigenvalue problem for a partial
differential operator reduces to an eigenvalue problem for an ordinary differential
operator. This case we call the \emph{axisymmetric case.}
\end{remark}
Consider therefore the perturbed metric $g_{\alpha\beta}(x^1;\upepsilon)$, the components of which are smooth functions of the coordinate $x^1$ and a small real parameter $\upepsilon$ which satisfies
\begin{equation}\label{unperturbed metric}
g_{\alpha\beta}(x^1;0)=\delta_{\alpha\beta}.
\end{equation}
For a given function $f:\mathbb{T}^3\rightarrow \mathbb{C}$ we denote its Fourier coefficients by
\begin{equation}\label{fourier coefficients}
\widehat{f}(m_1):=\frac{1}{2\pi}\int_{\mathbb{T}^3} e^{-im_1 x^1}f(x^1)dx^1,\ \ m_1\in \mathbb{Z},
\end{equation} and  we let
\begin{equation}\label{matrices h and k}
h_{\alpha\beta}(x^1):=\left. \frac{\partial g_{\alpha\beta}(x^1;\upepsilon)}{\partial \upepsilon}\right\vert_{\upepsilon=0},\quad \
k_{\alpha\beta}(x^1):=4\left. \frac{\partial^2 g_{\alpha\beta}(x^1;\upepsilon)}{\partial \upepsilon^2}\right\vert_{\upepsilon=0}.
\end{equation}
We choose to deal with the massless Dirac operator on half-densities denoted by $W_{1/2}$ which is the operator $\displaystyle W_{1/2}:=(\det g_{\kappa\lambda})^{1/4}W(\det g_{\mu\nu})^{-1/4}$, which therefore differs from the massless Dirac operator $W$ only by ``scalar'' factors on the left and on the right, as was explained in \cite{chervova2013spectral}.
It is known that these operators are equivalent, i.e. have the same spectrum, see \cite{chervova2013spectral, downes2013spectral}.

The axisymmetric massless Dirac operator on half-densities corresponding to the perturbed metric $g(x^1;\upepsilon)$ reads
\begin{align}
\nonumber
W_{1/2}(\upepsilon)=&-\frac{i}{2}
\left(\!\!
  \begin{array}{cc}
    \tensor{e}{_3^1} & \tensor{e}{_1^1}-i\tensor{e}{_2^1}\\
    \tensor{e}{_1^1}+i\tensor{e}{_2^1} & -\tensor{e}{_3^1} \\
  \end{array}\!\!
\right)\frac{d}{dx^1} -
\frac{i}{2}\frac{d}{dx^1}
\left(\!\!
  \begin{array}{cc}
    \tensor{e}{_3^1} & \tensor{e}{_1^1}-i\tensor{e}{_2^1}\\
    \tensor{e}{_1^1}+i\tensor{e}{_2^1} & -\tensor{e}{_3^1} \\
  \end{array}\!\!
\right)\label{massless dirac}\\
&+\frac{\delta_{jk}}{4\sqrt{\det g_{\alpha\beta}}}
\left(
\tensor{e}{^j_3}\left(\frac{d\tensor{e}{^k_2}}{dx^1}\right)
-\tensor{e}{^j_2}\left(\frac{d\tensor{e}{^k_3}}{dx^1}\right)
\right)I,
\end{align}
where ${e}{^i{}_j}$ and ${e}{_i{}^j}$ are the components of the perturbed coframe $e^j{}_\alpha(x^1;\upepsilon)$ and perturbed frame $e_j{}^\alpha(x^1;\upepsilon)$ respectively,
$I$ is identity matrix  is $2\times 2$ and
\begin{equation*}
\sqrt{\det g_{\alpha\beta}}=\frac{1}{\sqrt{\det g^{\alpha\beta}}}=
\det \tensor{e}{^j_\alpha}=\frac{1}{\det \tensor{e}{_j^\alpha}},
\end{equation*} see Section 7 of \cite{downes2013spectral}, where it was explained that the operator $W_{1/2}(\upepsilon)$ acts on 2-columns
$v=(v_1,\  v_2)^T$
of complex-valued half-densities. Our
Hilbert space
is the vector space of 2-columns of square integrable half-densities equipped with inner product
\begin{equation}\label{scalar product}
\langle v, w\rangle:=\int_M w^\ast v dx\ .
\end{equation}
The domain of the operator $W_{1/2}(\upepsilon)$ is $H^1(M;\mathbb{C}^2)$, which is the Sobolev space of 2-columns of half-densities that are square integrable together with their first partial derivatives.
\begin{remark}\label{remarkDiscreteSpectrum}
The operator $W_{1/2}(\upepsilon):H^1(M;\mathbb{C}^2)\rightarrow L^2(M;\mathbb{C}^2)$ is self-adjoint and its spectrum is discrete, with eigenvalues accumulating to $\pm \infty$, see \cite{downes2013spectral}.
\end{remark}
One of the main reasons we choose to work with the massless Dirac operator
on half-densities $W_{1/2}(\upepsilon)$ rather than with the massless Dirac operator $W(\upepsilon)$ is that we do not want our Hilbert space to depend on $\upepsilon$.
According to formulae (\ref{unperturbed metric}) and (\ref{matrices h and k}) we have that
\begin{equation}\label{metric-asymptotic expansion}
g_{\alpha\beta}(x^1;\upepsilon)=\delta_{\alpha\beta}+\upepsilon h_{\alpha\beta}(x^1)+\frac{\upepsilon^2}{4} k_{\alpha\beta}(x^1)+O(\upepsilon^3)
\end{equation}
and hence, using Taylor's formula for the function $\sqrt{1+z}$, the coframe is given by
\begin{equation}\label{coframe}
e^j{}_\alpha(x^1;\upepsilon)=\delta^j{}_\alpha+\frac{\upepsilon}{2}h^j{}_\alpha-\frac{\upepsilon^2}{8}
\left(h^2\right)^j{}_\alpha+\frac{\upepsilon^2}{8}k^j{}_\alpha+O(\upepsilon^3)
\end{equation}
and, using Taylor's formula for the function $(1+z)^{-1}$, the frame is given by
\begin{equation}\label{frame}
e_j{}^\alpha(x^1;\upepsilon)=\delta_j{}^\alpha
-\frac{\upepsilon}{2}h_j{}^\alpha+\frac{3\upepsilon^2}{8}\left(h^2\right)_j{}^\alpha - \frac{\upepsilon^2}{8}k_j{}^\alpha + O(\upepsilon^3).
\end{equation}
\begin{remark}Note that $h^2$  denotes the square of the perturbation matrix  $h$, i.e.  $\left(h^2\right)_j{}^\alpha=\tensor{h}{_j^\beta}\tensor{h}{_\beta^\alpha}$, where the summation is performed over the repeated index $\beta$.
\end{remark}
Note that for a given metric $g_{\alpha\beta}(x^1;\upepsilon)$ the coframe $\tensor{e}{^j_\alpha}(x^1;\upepsilon)$ and the frame $\tensor{e}{_j^\alpha}(x^1;\upepsilon)$ are not defined uniquely.
We can multiply the matrix functions $\tensor{e}{^j_\alpha}(x^1;\upepsilon)$ and $\tensor{e}{_j^\alpha}(x^1;\upepsilon)$ from the left by an arbitrary smooth $3\times 3$ special orthogonal matrix-function $O(x^1;\upepsilon)$ satisfying the condition $O(x^1;0)=I$. It is important to stress that this choice of the coframe does not affect the spectrum of the massless Dirac operator, see \cite{chervova2013spectral,downes2013spectral} for more details.

It is also known that the eigenfunctions of the operator
$W_{1/2}(\upepsilon)$ are infinitely smooth, see
\cite{chervova2013spectral,downes2013spectral,downes2016spectral}.
We get directly that for any integer eigenvalue  $\lambda\in \mathbb{Z}$, see Remark \ref{remarkDiscreteSpectrum}, the corresponding eigenvector of the massless Dirac operator $W_{1/2}(0)$ is
\begin{align}\label{eigenvector}
v_\lambda(x^1)&=\frac{1}{2\sqrt{\pi}}\left(
      \begin{array}{c}
        1 \\
        1 \\
      \end{array}
    \right)e^{i\lambda x^1}.
\end{align}
\begin{remark}\label{remarkEigenvector}
As each eigenvalue of massless Dirac operator in the axisymmetric case has multiplicity two, see \cite{downes2013spectral} and we also have that the vector
\begin{align}\label{eigenvector Cv}
w_\lambda(x^1)&=C(v_\lambda(x^1))  =\frac{1}{2\sqrt{\pi}}\left(
      \begin{array}{r}
        -1 \\
        1 \\
      \end{array}
    \right)e^{-i\lambda x^1}
\end{align}
is also an eigenvector of the massless Dirac operator corresponding to eigenvalue $\lambda\in\mathbb{Z}$, where $C$ is the operator of charge conjugation (\ref{eigenvector Cv2}).
\end{remark}
Let $W_{1/2}(\upepsilon)$ therefore be the massless Dirac operator on half-densities \eqref{massless dirac} corresponding to the perturbed metric $g_{\alpha\beta}(x;\upepsilon)$ and let
\begin{equation}\label{asymptotic expansion of the operator}
W_{1/2}(\upepsilon)=W_{1/2}^{(0)}+\upepsilon W_{1/2}^{(1)}+\upepsilon^2 W_{1/2}^{(2)}+\cdots
\end{equation}
be the asymptotic expansion of the perturbed massless Dirac operator in powers of the small
parameter $\upepsilon$. The operator $W_{1/2}^{(0)}=W_{1/2}(0)$ is the unperturbed massless Dirac operator on half-densities \eqref{massless dirac}. We denote by $\lambda^{(0)}$  the eigenvalue of this operator and by $v^{(0)}$ the corresponding eigenvector.

The perturbation of an isolated eigenvalue of finite
multiplicity of a bounded operator  was described in \cite{rellich1969perturbation} and that procedure can be applied in our case with some additional conditions.

The operators $W_{1/2}^{(k)}$, $k=0,1,2,\ldots$ are formally self-adjoint first order
differential operators which also commute with the antilinear operator of charge conjugation \eqref{eigenvector Cv2}.
We need to solve the eigenvalue problem
$$
W_{1/2}(\upepsilon)v(\upepsilon)=\lambda(\upepsilon) v(\upepsilon).
$$
We seek the eigenvalue and eigenfunction of the perturbed operator $W_{1/2}(\upepsilon)$ in the
form of asymptotic expansions
\begin{align}\label{asimptotski razvoj lambde}
\lambda(\upepsilon)&=\lambda^{(0)}+\upepsilon \lambda^{(1)}+\upepsilon^2\lambda^{(2)}+\cdots,\\
\label{asimptotski razvoj v}
v(\upepsilon)&=v^{(0)}+\upepsilon v^{(1)}+\upepsilon^2v^{(2)}+\cdots.
\end{align}
In general, these asymptotic series do \emph{not} converge and the formal procedure for this perturbation process was described in Section 4 in \cite{downes2013spectral}.
However, the asymptotic expansion can be justified, as was described in Section 5 in \cite{downes2013spectral}, by showing that by taking a finite number of terms one gets the expected estimate for the remainder. This argument holds for any simple eigenvalue.

\section{Main results}
Let $n\in\mathbb{N}$ be a positive eigenvalue of the unperturbed massless Dirac operator. Throughout the rest of this paper we denote by  $\lambda_{+n}(\upepsilon)$ the asymptotic expansion of the eigenvalue $n$ and by $\lambda_{-n}(\upepsilon)$ the asymptotic expansion of the eigenvalue $-n$. By $\lambda^{(i)}_{+n}$ and $\lambda^{(i)}_{-n}$, $i=1,2,3,\ldots$, we denote their respective asymptotic coefficients.

The main result of this paper is the following
\begin{theorem}\label{main theorem}
Under an arbitrary perturbation of the metric  (\ref{metric-asymptotic expansion}), the asymptotic expansion of the eigenvalues $n$ and $-n$ are
\begin{equation}\label{asymptotic expansion of eigenvalues}
\lambda_{\pm n}(\upepsilon)=\pm n +\lambda^{(1)}_{\pm n}\upepsilon+\lambda^{(2)}_{\pm n}\upepsilon^2+O(\upepsilon^3)\ \ \textrm{as}\ \ \upepsilon \rightarrow 0,
\end{equation}
where the constants $\lambda_{+n}^{(1)}$, $\lambda_{-n}^{(1)}$, $\lambda_{+n}^{(2)}$ and $\lambda_{-n}^{(2)}$ appearing in \eqref{asymptotic expansion of eigenvalues} are given by
\begin{equation}\label{lambda1NEW}
\lambda^{(1)}_{\pm n}=\mp \frac{n}{2}\tensor{\widehat{h}}{_{11}}(0),
\end{equation}
\begin{align}
\label{lambda2NEW}
\lambda^{(2)}_{\pm n}&=
\pm \frac{3n}{8} \widehat{(h^2)}{}_{11}(0)
\mp \frac{n}{8}\tensor{\widehat{k}}{_{11}}(0)
-
\frac{i}{16}\varepsilon_{\beta\gamma1}\sum_{m\in \mathbb{Z}\setminus \{0\}}
m\overline{\widehat{h}_{\alpha\beta}(m)}\widehat{h}_{\alpha\gamma}(m)
\\
& -\frac{1}{16}\sum_{m\in \mathbb{Z}\setminus \{\pm n\}}
\frac{1}{m\mp n}
(m\pm n)^2\ \widehat{h}_{11}(m\mp n)
\overline{\widehat{h}_{11}(m\mp n)},
\nonumber
\\
 \nonumber
&
-\frac{1}{16}\sum_{m\in \mathbb{Z}\setminus \{\pm n\}}
(m\mp n)\
\left(\widehat{h}_{31}(m\pm n)
+i\widehat{h}_{21}(m\pm n)\right) \\
\nonumber
& ~ \hspace{4cm}
\left(\overline{\widehat{h}_{31}(m\pm n)}
-i\overline{\widehat{h}_{21}(m\pm n)}\right),
\end{align}
\end{theorem}
\begin{remark}
Note that, as usual,  $\varepsilon_{\alpha\beta\gamma}$ denotes the totally antisymmetric quantity, $\varepsilon_{123}:=+1$, while overline stands for complex conjugation.
\end{remark}
\begin{remark}
Note that \eqref{lambda1NEW} directly implies that spectral asymmetry cannot be achieved in the linear term.
\end{remark}

%
%

\begin{remark}
Formula (\ref{lambda2NEW}) implies that
\begin{align*}
\nonumber
&\lambda^{(2)}_{+n}+\lambda^{(2)}_{-n}=
-
\frac{i}{8}\varepsilon_{\beta\gamma1}\sum_{m\in \mathbb{Z}\setminus \{0\}}
m\overline{\widehat{h}_{\alpha\beta}(m)}\widehat{h}_{\alpha\gamma}(m)\\
\nonumber
&-\frac{1}{16}\sum_{m\in \mathbb{Z}\setminus \{n\}} \frac{1}{m-n}(m+n)^2\
\widehat{h}_{11}(m-n)\overline{\widehat{h}_{11}(m-n)}\\
\nonumber
&-\frac{1}{16}\sum_{m\in \mathbb{Z}\setminus \{n\}} (m-n)
\left(\widehat{h}_{31}(m+n)+i\widehat{h}_{21}(m+n)\right)
\left(\overline{\widehat{h}_{31}(m+n)}
-i\overline{\widehat{h}_{21}(m+n)}\right)\\
&-\frac{1}{16}\sum_{m\in \mathbb{Z}\setminus \{-n\}}
\frac{1}{m+n}
(m-n)^2\
\widehat{h}_{11}(m+n)\overline{\widehat{h}_{11}(m+n)}\\
&-\frac{1}{16}\sum_{m\in \mathbb{Z}\setminus \{-n\}}(m+n)
\left(\widehat{h}_{31}(m-n)+i\widehat{h}_{21}(m-n)\right)
\left(\overline{\widehat{h}_{31}(m-n)}
-i\overline{\widehat{h}_{21}(m-n)}\right).
\end{align*}
We conclude that one way to create spectral asymmetry is to choose the perturbation matrix $h_{\alpha\beta}(x^1)$ such that
$\widehat{h}_{11}(m\pm n) = \widehat{h}_{21}(m\pm n)
=\widehat{h}_{31}(m\pm n)=0$ for all $m\in\mathbb{Z}$ such that the term $$\displaystyle \varepsilon_{\beta\gamma1}\sum_{m\in \mathbb{Z}\setminus \{0\}}
m\overline{\widehat{h}_{\alpha\beta}(m)}\widehat{h}_{\alpha\gamma}(m)\neq 0.$$
\end{remark}

\begin{remark}
The length of the unit $x^1$ circle is  $\displaystyle 2\pi\left(1+\frac{1}{2}\widehat{h}_{11}(0)\upepsilon\right)+O(\upepsilon^2).$ Therefore the coefficients $\lambda^{(1)}_{+n}(\upepsilon)$ and $\lambda^{(1)}_{-n}(\upepsilon)$ are determined by the change of the length of the unit $x^1$ circle.
\end{remark}

\begin{proof}[Proof of Theorem \ref{main theorem}.]
In order to prove Theorem  \ref{main theorem} we need to write down explicitly the differential operators
$W^{(1)}$ and $W^{(2)}$ appearing in the asymptotic expansion of the perturbed massless Dirac operator on half-densities
$$
W_{1/2}(\upepsilon)=W^{(0)}+\upepsilon W^{(1)}+\upepsilon^2W^{(2)}+O(\upepsilon^3).
$$
Substituting formulae (\ref{coframe}) and (\ref{frame}) into (\ref{massless dirac}), we get that
\begin{align}
\label{A1}
W^{(1)}&=\frac{i}{4}
\left(\begin{array}{cc}
                      \tensor{h}{_3^1} & \tensor{h}{_1^1}-i\tensor{h}{_2^1}\\
                      \tensor{h}{_1^1}+i\tensor{h}{_2^1} & -\tensor{h}{_3^1}
                    \end{array}
\right)\frac{d}{dx^1}\\
\nonumber
&+
\frac{i}{4}\frac{d}{dx^1}
\left(\begin{array}{cc}
                      \tensor{h}{_3^1} & \tensor{h}{_1^1}-i\tensor{h}{_2^1}\\
                      \tensor{h}{_1^1}+i\tensor{h}{_2^1} & -\tensor{h}{_3^1}
                    \end{array}
\right),
\end{align}

\begin{align}\nonumber
W^{(2)}&=-\frac{3i}{16}
\left(\begin{array}{cc}
                      \left(h^2\right)_3{}^1 & \left(h^2\right)_1{}^1-i \left(h^2\right)_2{}^1\\
                      \left(h^2\right)_1{}^1+i \left(h^2\right)_2{}^1 & -\left(h^2\right)_3{}^1
                    \end{array}
\right)\frac{d}{dx^1}\\
\nonumber
&-\frac{3i}{16}\frac{d}{dx^1}
\left(\begin{array}{cc}
                      \left(h^2\right)_3{}^1 & \left(h^2\right)_1{}^1-i \left(h^2\right)_2{}^1\\
                      \left(h^2\right)_1{}^1+i \left(h^2\right)_2{}^1 & -\left(h^2\right)_3{}^1
                    \end{array}
\right)
-\frac{1}{16}\varepsilon_{\beta\gamma1}h_{\alpha\beta}\frac{dh_{\alpha\gamma}}{dx^1}I.
\\
&+
\frac{i}{16}
\left(\begin{array}{cc}
                      \tensor{k}{_3^1} & \tensor{k}{_1^1}-i\tensor{k}{_2^1}\\
                      \tensor{k}{_1^1}+i\tensor{k}{_2^1} & -\tensor{k}{_3^1}
                    \end{array}
\right)\frac{d}{dx^1}
+
\frac{i}{16}\frac{d}{dx^1}
\left(\begin{array}{cc}
                      \tensor{k}{_3^1} & \tensor{k}{_1^1}-i\tensor{k}{_2^1}\\
                      \tensor{k}{_1^1}+i\tensor{k}{_2^1} & -\tensor{k}{_3^1}
                    \end{array}
\right)
\label{A2}
\end{align}
We calculate the coefficients $\lambda^{(1)}$ and $\lambda^{(2)}$ from equation \eqref{asimptotski razvoj lambde} using
\begin{equation}\label{formula za lambda 1}
\lambda^{(1)}=\langle W^{(1)} v^{(0)},v^{(0)}  \rangle,
\end{equation}
\begin{equation}\label{formula za lambda 2}
\lambda^{(2)}=\langle W^{(2)}v^{(0)},v^{(0)}\rangle
-\langle(W^{(1)}-\lambda^{(1)})Q(W^{(1)}-\lambda^{(1)})v^{(0)},v^{(0)}\rangle,
\end{equation}
see Appendix \ref{appendix explicit formulae} for details, where $v^{(0)}$ is the eigenvector corresponding to the eigenvalue  $\lambda^{(0)}$. Note that $Q$ denotes the \emph{pseudoinverse operator} (see Appendix \ref{Section Pseudoinverse operator construction})
of the unperturbed operator $W_{1/2}(0)$ corresponding to the integer eigenvalue $\lambda$. It is given by the formula
\begin{align}
\label{pseudoinverse operator}
Q=\frac{1}{4\pi}
\sum_{m\in\mathbb{Z}\setminus \{\lambda\}}\frac{1}{m-\lambda}
&\left[
e^{imx^1}
\left(
  \begin{array}{cc}
    1 & 1 \\
    1 & 1 \\
  \end{array}
\right)\int_{0}^{2\pi}e^{-imy^1}(\ \cdot\ )dy^1\right.\\
&\left.+
e^{-imx^1}
\left(
  \begin{array}{cc}
    1 & -1 \\
    -1 & 1 \\
  \end{array}
\right)\int_{0}^{2\pi}e^{imy^1}(\ \cdot\ )dy^1
\right],\nonumber
\end{align}
see e.g. \cite{downes2013spectral,rellich1969perturbation} and Appendix \ref{Section Pseudoinverse operator construction} for the general details on the pseudoinverse operator.\\
Using formula (\ref{fourier coefficients}) for the Fourier coefficients, substituting formulae  (\ref{eigenvector}), (\ref{A1}), (\ref{A2}) and (\ref{pseudoinverse operator}) into (\ref{formula za lambda 1}) and (\ref{formula za lambda 2}),  we obtain equations (\ref{lambda1NEW}) and (\ref{lambda2NEW}), which describe the linear and quadratic terms in the asymptotic expansion of eigenvalues $\pm n$.  Note that we get the same formulae if we use the eigenvector (\ref{eigenvector Cv}).
\end{proof}

\begin{remark}
Appendix \ref{AppendixSecondTermMasslessDirac} contains very detailed explicit derivations of equations \eqref{lambda1NEW} and \eqref{lambda2NEW}.
\end{remark}
\section{Explicit examples of spectral asymmetry}
In this section we also give two explicit families of the perturbed Euclidean metric for which we obtain spectral asymmetry of the massless Dirac operator in the axisymmetric case.
\begin{example}\label{example1 asymmetry} If we consider
\[
 h_{\alpha\beta}(x^1)=2\left(
     \begin{array}{ccc}
       0 & 0 & 0 \\
       0 & \cos x^1 & \sin x^1 \\
       0 & \sin x^1 & -\cos x^1 \\
     \end{array}
   \right), \quad
k_{\alpha\beta}(x^1)=\left(
    \begin{array}{ccc}
      \sin x^1 & \cos x^1 & 0 \\
      \cos x^1 & 0 & 0 \\
      0 & 0 & 0 \\
    \end{array}
  \right),
\]
then by using formulae (\ref{A1}) and (\ref{A2}), we get that
the corresponding perturbed massless Dirac operator is given by
\begin{align*}
W(\upepsilon)=&-i
\left(
  \begin{array}{cc}
    0 & 1 \\
    1 & 0 \\
  \end{array}
\right)\frac{d}{dx^1}+
\frac{\upepsilon^2}{16}
\left(
  \begin{array}{cc}
    0 & \cos x^1+i\sin x^1 \\
    -\cos x^1+i\sin x^1 & 0 \\
  \end{array}
\right)\frac{d}{dx^1}\\
&+
\frac{\upepsilon^2}{16}\frac{d}{dx^1}
\left(
  \begin{array}{cc}
     0 & \cos x^1+i\sin x^1 \\
    -\cos x^1+i\sin x^1 & 0 \\
  \end{array}
\right)-\frac{\upepsilon^2}{2} I.
\end{align*}
Using formula (\ref{pseudoinverse operator}) for the pseudoinverse operator
   in the case of eigenvalues $n=\pm 1$,  as well as  (\ref{scalar product}), (\ref{formula za lambda 1}) and (\ref{formula za lambda 2}),  we get that
\begin{align}\label{first example eigenvalue one}
\lambda_{+1}(\upepsilon)&=1-\frac{1}{2}\upepsilon^2+O(\upepsilon^3),\\
\label{first example eigenvalue minus one}
\lambda_{-1}(\upepsilon)&=-1-\frac{1}{2}\upepsilon^2+O(\upepsilon^3).
\end{align}
Applying the Fourier transform (\ref{fourier coefficients}) to matrices $h$ and $k$ gives us
\begin{equation}\label{example1 fourier h}
\widehat{h}_{\alpha\beta}(m)=\left\{
\begin{array}{ccl}
  \left(
    \begin{array}{ccc}
      0 & 0 & 0 \\
      0 & 1 & -i \\
      0 & -i & -1 \\
    \end{array}
  \right)
   &\textrm{for}& m=1, \\ \\
  0 & \textrm{for} &m=2,3,\ldots,
\end{array}
\right.
\end{equation}
\begin{equation}\label{example1 fourier k}
\widehat{k}_{\alpha\beta}(m)=\left\{
\begin{array}{ccl}
  \left(
    \begin{array}{ccc}
      -i/2 & 1/2 & 0 \\
      1/2 & 0 & 0 \\
      0 & 0 & 0 \\
    \end{array}
  \right)
   &\textrm{for}& m=1, \\ \\
  0 & \textrm{for} &m=2,3,\ldots,
\end{array}
\right.
\end{equation}
\begin{equation}\label{example1 fourier hsquared}
\widehat{(h^2)}_{\alpha\beta}(m)=\left\{
\begin{array}{ccl}
  \left(
    \begin{array}{ccc}
      0 & 0 & 0 \\
      0 & 4 & 0 \\
      0 & 0 & 4 \\
    \end{array}
  \right)
   &\textrm{for}& m=0, \\ \\
  0 & \textrm{for} &m=1,2,3,\ldots
\end{array}
\right.
\end{equation}
Substituting (\ref{example1 fourier h}), (\ref{example1 fourier k}) and (\ref{example1 fourier hsquared}) into (\ref{lambda1NEW}) and (\ref{lambda2NEW}), we get that $\lambda^{(1)}_{\pm 1}=0$ and $ \lambda^{(2)}_{\pm 1}=-\frac12$,  which is in accordance with (\ref{first example eigenvalue one}) and (\ref{first example eigenvalue minus one}).
\end{example}

\begin{example}\label{example2 asymmetry} If we consider
\[
h_{\alpha\beta}(x^1)=\left(
     \begin{array}{ccc}
       1 & \cos x^1 & \sin x^1 \\
       \cos x^1 & \cos x^1 & \sin x^1 \\
       \sin x^1 & \sin x^1 & -\cos x^1 \\
     \end{array}
   \right), \quad
k_{\alpha\beta}(x^1)=\left(
    \begin{array}{ccc}
      \sin x^1 & \cos x^1 & 0 \\
      \cos x^1 & -\sin x^1 & 0 \\
      0 & 0 & 0 \\
    \end{array}
  \right),
\]
then using formulae (\ref{A1}) and (\ref{A2}), we get that
the corresponding perturbed massless Dirac operator is given by
\begin{align*}
W(\upepsilon)&=-i
\left(
  \begin{array}{cc}
    0 & 1 \\
    1 & 0 \\
  \end{array}
\right)\frac{d}{dx^1}
+\frac{\upepsilon}{4}
\left(
  \begin{array}{cc}
    i\sin x^1 & i+\cos x^1 \\
    i-\cos x^1 & -i\sin x^1 \\
  \end{array}
\right)
\frac{d}{dx^1} \\
&
+\frac{\upepsilon}{4}\frac{d}{dx^1}
\left(
  \begin{array}{cc}
    i\sin x^1 & i+\cos x^1 \\
    i-\cos x^1 & -i\sin x^1 \\
  \end{array}
\right) -\frac{3\upepsilon^2}{16}I
\\
&+\frac{\upepsilon^2}{16}
\left(
  \begin{array}{cc}
    -3i\sin x^1 & -6i-3-2\cos x^1+i\sin x^1 \\
    -6i+3+2\cos x^1+i\sin x^1 & 3i\sin x^1 \\
  \end{array}
\right)
\frac{d}{dx^1}\\
\nonumber
&+\frac{\upepsilon^2}{16}\frac{d}{dx^1}
\left(
  \begin{array}{cc}
    -3i\sin x^1 & -6i-3-2\cos x^1+i\sin x^1 \\
    -6i+3+2\cos x^1+i\sin x^1 & 3i\sin x^1 \\
  \end{array}
\right).
\end{align*}
Using formula (\ref{pseudoinverse operator}) for the pseudoinverse operator in the case of eigenvalues  $n=\pm 1$, as well as (\ref{scalar product}), (\ref{formula za lambda 1}) and (\ref{formula za lambda 2}),  we get that
\begin{align}\label{second example eigenvalue one}
\lambda_{+1}(\upepsilon)&=1-\frac12\upepsilon+\frac34\upepsilon^2+O(\upepsilon^3),\\
\label{second example eigenvalue minus one}
\lambda_{-1}(\upepsilon)&=-1+\frac12\upepsilon-\upepsilon^2+O(\upepsilon^3).
\end{align}
Application of the Fourier transform (\ref{fourier coefficients}) gives us
\begin{equation}\label{example2 fourier h}
\widehat{h}_{\alpha\beta}(m)=\left\{
\begin{array}{ccl}
\left(
    \begin{array}{ccc}
      1 & 0 & 0 \\
      0 & 0 & 0 \\
      0 & 0 & 0 \\
    \end{array}
  \right)
   &\textrm{for}& m=0,\\ \\
  \left(
    \begin{array}{ccc}
      0 & 1/2 & -i/2\\
      1/2 & 1/2 & -i/2 \\
      -i/2 & -i/2 & -1/2 \\
    \end{array}
  \right)
   &\textrm{for}& m=1, \\ \\
  0 & \textrm{for} &m=2,3,\ldots,
\end{array}
\right.
\end{equation}
\begin{equation}\label{example2 fourier k}
\widehat{k}_{\alpha\beta}(m)=\left\{
\begin{array}{ccl}
  \left(
    \begin{array}{ccc}
      -i/2 & 1/2 & 0 \\
      1/2 & i/2 & 0 \\
      0 & 0 & 0 \\
    \end{array}
  \right)
   &\textrm{for}& m=1, \\ \\
  0 & \textrm{for} &m=2,3,\ldots,
\end{array}
\right.
\end{equation}
\begin{equation}\label{example2 fourier hsquared}
\widehat{(h^2)}_{\alpha\beta}(m)=\left\{
\begin{array}{ccl}
\left(
    \begin{array}{ccc}
      2 & 1 & 0 \\
      1 & 3/2 & 0 \\
      0 & 0 & 3/2 \\
    \end{array}
  \right)
   &\textrm{for}& m=0,\\ \\
  \left(
    \begin{array}{ccc}
      0 & 1/2 & -i/2 \\
      1/2 & 0 & 0 \\
      -i/2 & 0 & 0 \\
    \end{array}
  \right)
   &\textrm{for}& m=1, \\ \\
   \left(
    \begin{array}{ccc}
      0 & 0 & 0 \\
      0 & 1/2 & -i/4 \\
      0 & -i/4 & -1/4 \\
    \end{array}
  \right)
   &\textrm{for}& m=2,\\ \\
  0 & \textrm{for} &m=3,4,\ldots.
\end{array}
\right.
\end{equation}
Substituting (\ref{example2 fourier h}), (\ref{example2 fourier k}) and (\ref{example2 fourier hsquared}) into (\ref{lambda1NEW}) and (\ref{lambda2NEW}),
we get that $\displaystyle \lambda^{(1)}_{\pm 1}=\mp\frac12$, $\displaystyle \lambda^{(2)}_{+1}=\frac34$ and $\displaystyle \lambda^{(2)}_{-1}=-1$  which is in accordance with (\ref{second example eigenvalue one}) and (\ref{second example eigenvalue minus one}).
\end{example}

\section{Conclusions and future goals}
In this paper we obtained spectral asymmetry of the massless Dirac operator in the axisymmetric case for an arbitrary eigenvalue different from zero. We also developed general methods and performed calculations which we can apply in order to obtain spectral asymmetry not only in the axisymmetric case but also in the general case \eqref{dirac operator with euclidean metric}, which we were not able to do at this time, but hope to achieve in the near future.

The eigenvalues with smallest modulus of the massless Dirac operator on the unit torus $\mathbb{T}^3$ equipped with Euclidean metric and standard spin structure and the unit sphere $\mathbb{S}^3$ equipped with Riemannian metric were considered in \cite{downes2013spectral} and \cite{fang2016spectral}, for which the asymptotic formulae were derived. In this paper  we worked on the unit torus $\mathbb{T}^3$ and for the axisymmetric case we derived the asymptotic formulae for the eigenvalues $\pm n$.

Our future goal is to consider eigenvalues, which are not with smallest modulus, of the massless Dirac operator when the metric not only depends on one spatial coordinate and to derive asymptotic formulae for these eigenvalues in order to obtain spectral asymmetry.

\section*{Appendices}
\appendix

\section{Pseudoinverse operator construction}
\label{Section Pseudoinverse operator construction}
Let $v^{(0)}$ be a normalized eigenvector of the operator $A$ corresponding to the eigenvalue $\lambda^{(0)}$. Then the vector $C(v^{(0)})$ is also a normalized eigenvector corresponding to the eigenvalue $\lambda^{(0)}$, see Remark \ref{remarkEigenvector}.
Consider the problem
\begin{equation}\label{temp1}
\left(A-\lambda^{(0)}\right)v=f,
\end{equation}
for a given function $f\in L^2(M;\mathbb{C}^2)$ where we need to find the function $v\in H^1(M;\mathbb{C}^2)$. Suppose that the function $f$ satisfies the conditions
$$
\langle f, v^{(0)}\rangle=\langle f, C(v^{(0)})\rangle=0,
$$
where  $C$ is the charge conjugation operator \eqref{eigenvector Cv2}. The  problem \eqref{temp1} can be resolved for the function $v$ and the uniqueness of this function is achieved by the conditions
$$
\langle v, v^{(0)}\rangle=\langle v, C(v^{(0)})\rangle=0.
$$
We define the operator $Q$ as
$Q:f\mapsto v.$
The operator $Q$ is a bounded linear operator acting on the orthogonal complement of the eigenspace of the
operator $A$ corresponding to the eigenvalue $\lambda^{(0)}$. Also, the bounded linear operator $Q$ is self-adjoint and commutes with the antilinear operator of charge conjugation (\ref{eigenvector Cv2}).
We can extend the acting of the pseudoinverse operator $Q$  to the whole Hilbert
space $L^2(M;\mathbb{C}^2)$ in accordance with
$
Qv^{(0)}=QC(v^{(0)})=0.
$

Now we present the construction of the pseudoinverse operator $Q$  itself,  largely following the exposition of \cite{rellich1969perturbation}.
Let $\lambda$ be an eigenvalue of multiplicity $k$ of a  Hermitian operator $A$.
Then the homogenous equation
$
(A-\lambda^{(0)})v=0
$
has $k$ linearly independent solutions $\phi^{(1)},\phi^{(2)},\ldots,\phi^{(k)}$ for which we can assume  orthonormality, i.e.
\begin{equation*}
\langle\phi^{(i)},\phi^{(j)}\rangle=\delta_{ij},\ \ (i,j=1,2,\ldots,k).
\end{equation*}
The operator $A-\lambda^{(0)}$ has no inverse, but there is a unique bounded Hermitian operator $Q$ such that
 $Q\phi^{(i)}=0,$ $ (i=1,2,\ldots, k )$ and
 $$\displaystyle Q(A-\lambda^{(0)})u=u-\sum_{i=1}^{k}\langle\phi^{(i)},u\rangle\phi^{(i)}.$$
Define the projector operator $P$ by
\begin{equation*}\label{ProjectionOperatorP}
Pu :=\sum_{i=1}^{k}\langle\phi^{(i)},u\rangle\phi^{(i)}.
\end{equation*}
Then the properties of the above operator $Q$  can be written as $QP=0$ and
 $Q(A-\lambda^{(0)})=I-P$. The operator $Q$ is called the pseudoinverse operator of  the operator $A-\lambda.$
We can complete the eigenvectors $\phi^{(1)},\phi^{(2)},\ldots,\phi^{(k)}$ with eigenvectors $\phi^{(k+1)}, \phi^{(k+2)}, \ldots, \phi^{(n)}$, which correspond to eigenvalues $\lambda_{k+1},\lambda_{k+2},\ldots,\lambda_{n}$, respectively, such that
$$\langle\phi^{(i)},\phi^{(j)}\rangle=\delta_{ij},\ \ (i,j=1,2,\ldots,n).$$
Expanding $v$ and $f$ as
\begin{equation*}
v=\sum_{i=1}^{n}\langle\phi^{(i)},v\rangle\phi^{(i)},\ \  f=\sum_{i=1}^{n}\langle\phi^{(i)},f\rangle\phi^{(i)},
\end{equation*}
from  equation (\ref{temp1}) we get that
\begin{align*}
\langle\phi^{(i)},f\rangle&=0,\ \ (i=1,2,\ldots, k),\\
\langle\phi^{(i)},f\rangle&=(\lambda_i-\lambda)\langle\phi^{(i)},v\rangle,\ \ (i=k+1,k+2,\ldots, n).
\end{align*}
If the equation for $v$ has a solution, it is necessary for $f$
to be orthogonal to all solutions of the homogeneous equation. Hence we set
\begin{equation*}
v=\sum_{i=1}^{k}v_i\phi^{(i)}+\sum_{\lambda_i\neq \lambda}\frac{\langle\phi^{(i)},f\rangle}{\lambda_i- \lambda}\phi^{(i)},
\end{equation*}
where $v_i$ are arbitrary constants. The vector $v$ defines the complete solution of the equation.

Let the operator $P$ be the projector operator \ref{ProjectionOperatorP} into the space spanned by  the vectors $\phi^{(1)},\phi^{(2)},\ldots,\phi^{(k)}$ and $P_i$ the projector into the one-dimensional space spanned by $\phi^{(i)},\ (i=k+1,k+2,\ldots,n).$
\begin{definition}
The operator
\begin{equation*}
Q=\sum_{\lambda_i\neq \lambda^{(0)}}\frac{P_i}{\lambda_i-\lambda^{(0)}}
\end{equation*}
is the \emph{pseudoinverse operator} of the operator $A-\lambda^{(0)}$.
\end{definition}

\section{Explicit formulae for the asymptotic coefficients}\label{appendix explicit formulae}
Now we will derive the explicit formulae for the coefficients $\lambda^{(1)}$ and $\lambda^{(2)}$ in the asymptotic expansion \eqref{asimptotski razvoj lambde}.
Consider the perturbed eigenvalue problem
$$
W_{1/2}(\upepsilon)v(\upepsilon)=\lambda(\upepsilon) v(\upepsilon).
$$
Using formulae \eqref{asymptotic expansion of the operator}, \eqref{asimptotski razvoj lambde} and \eqref{asimptotski razvoj v}, we get that
\begin{align*}(W_{1/2}^{(0)}&+\upepsilon W_{1/2}^{(1)}+\upepsilon^2 W_{1/2}^{(2)}+\cdots)
(v^{(0)}+\upepsilon v^{(1)}+\upepsilon^2 v^{(2)}+\cdots)\\&=
(\lambda^{(0)}+\upepsilon \lambda^{(1)}+\upepsilon^2 \lambda^{(2)}+\cdots)
(v^{(0)}+\upepsilon v^{(1)}+\upepsilon^2 v^{(2)}+\cdots).
\end{align*}
By grouping together the elements not containing $\upepsilon$, we get that
$$W_{1/2}^{(0)}v^{(0)}=\lambda^{(0)}v^{(0)},$$
which is the unperturbed eigenvalue problem.
By grouping together the elements containing $\upepsilon$, we get that
$$
W_{1/2}^{(0)}v^{(1)}+W_{1/2}^{(1)}v^{(0)}=\lambda^{(0)}v^{(1)}+\lambda^{(1)}v^{(0)}
$$
and hence
$
(W_{1/2}^{(0)}-\lambda^{(0)})v^{(1)}=(\lambda^{(1)}-W_{1/2}^{(1)})v^{(0)}
$,
i.e. $v^{(1)}=Q(\lambda^{(1)}-W_{1/2}^{(1)})v^{(0)}$ where $Q$ is the pseudoinverse operator of the operator $W_{1/2}^{(0)}-\lambda^{(0)}.$
We denote by
\begin{equation}\label{tempf1}
f^{(1)}=(\lambda^{(1)}-W_{1/2}^{(1)})v^{(0)}.
\end{equation}
By grouping together the elements containing $\upepsilon^2$, we get that
$$
W_{1/2}^{(0)}v^{(2)}+W_{1/2}^{(1)}v^{(1)}+W_{1/2}^{(2)}v^{(0)}=
\lambda^{(0)}v^{(2)}+\lambda^{(1)}v^{(1)}+\lambda^{(2)}v^{(0)}
$$
and hence
$$
(W_{1/2}^{(0)}-\lambda^{(0)})v^{(2)}=(\lambda^{(2)}-W_{1/2}^{(2)})v^{(0)}+(\lambda^{(1)}-W_{1/2}^{(1)})v^{(1)}.
$$
We denote by
\begin{align}\nonumber
f^{(2)}&=(\lambda^{(2)}-W_{1/2}^{(2)})v^{(0)}+(\lambda^{(1)}-W_{1/2}^{(1)})v^{(1)}\\
\label{tempf2}
&
=(\lambda^{(2)}-W_{1/2}^{(2)})v^{(0)}+(\lambda^{(1)}-W_{1/2}^{(1)})Q(\lambda^{(1)}-W_{1/2}^{(1)})v^{(0)}
.
\end{align}
Continuing this process, the vectors $f^{(k)}$ and the coefficients $\lambda^{(k)}$ are obtained  from the conditions
\begin{equation}\label{tempCondition}
\langle f^{(k)},v^{(0)}\rangle=0,
\end{equation}
\begin{equation}\label{tempfk}
\langle f^{(k)},C(v^{(0)})\rangle=0.
\end{equation}
\begin{remark}
The eigenvalues have even multiplicity, so the condition \eqref{tempfk} is an additional condition which need to be satisfied. This is the part where our perturbation process differs from the standard perturbation process for single eigenvalues.
\end{remark}
The components $v^{(k)}$ are given by
$$
v^{(k)}=Qf^{(k)},
$$
where $Q$ is the pseudoinverse operator of the operator $W_{1/2}^{(0)}-\lambda^{(0)}.$
Substituting formulae \eqref{tempf1}, \eqref{tempf2} into formulae \eqref{tempCondition}, \eqref{tempfk} we obtain formulae \eqref{formula za lambda 1}, \eqref{formula za lambda 2}.

\section{Detailed calculations of the asymptotic coefficients}
\label{AppendixSecondTermMasslessDirac}
In this appendix we provide detailed calculations of the formulae for the  asymptotic coefficients \eqref{lambda1NEW} and \eqref{lambda2NEW} from Theorem~\ref{main theorem} which correspond to the eigenvalues $\pm n$. We will use the perturbation theory which  is described in Section~\ref{SectionPerturbation} and the explicitly derived formulae for the asymptotic coefficients given by \eqref{formula za lambda 1}, \eqref{formula za lambda 2}. We will also use the concept of the pseudoinverse of the massless Dirac operator  whose construction is given in Appendix~\ref{Section Pseudoinverse operator construction}.
\subsection{The calculations of the  $\lambda^{(1)}_{\pm n}$ coefficients}
We first want to derive equation \eqref{lambda1NEW} from Theorem~\ref{main theorem}.
Using the formula for the eigenvector \eqref{eigenvector} corresponding to the eigenvalues $n$ and $-n$ respectively, as well as  formula \eqref{A1}  for the differential operator $W_{1/2}^{(1)}$, integrating by parts, we get that the equation \eqref{formula za lambda 1} for the eigenvalues  $\pm n$ becomes
\begin{align}
\label{lambda1PlusNEW}
\lambda^{(1)}_{\pm n}&=\langle W_{1/2}^{(1)} v^{(0)},v^{(0)}  \rangle
=\int_{0}^{2\pi}[v^{(0)}]^\ast W_{1/2}^{(1)}v^{(0)}dx^1\\
\nonumber
&=\frac{i}{4}\int_{0}^{2\pi}[v^{(0)}]^\ast
\left(\begin{array}{cc}
                      \tensor{h}{_3^1} & \tensor{h}{_1^1}-i\tensor{h}{_2^1}\\
                      \tensor{h}{_1^1}+i\tensor{h}{_2^1} & -\tensor{h}{_3^1}
                    \end{array}
\right)\frac{d}{dx^1}  v^{(0)}dx^1\\
\nonumber
&-\frac{i}{4}\int_{0}^{2\pi}\frac{d}{dx^1} [v^{(0)}]^\ast
\left(\begin{array}{cc}
                      \tensor{h}{_3^1} & \tensor{h}{_1^1}-i\tensor{h}{_2^1}\\
                      \tensor{h}{_1^1}+i\tensor{h}{_2^1} & -\tensor{h}{_3^1}
                    \end{array}
\right)v^{(0)}dx^1\\
&\nonumber
=\mp \frac{n}{4\pi}\int_0^{2\pi}\tensor{h}{_1^1}(x^1)dx^1
=\mp \frac{n}{2}\tensor{\widehat{h}}{_{11}}(0),
\end{align}

where $\widehat{h}(m)$ denotes the Fourier coefficient \eqref{fourier coefficients}  for the function $h(x^1)$ .

\subsection{The calculations of the  $\lambda^{(2)}_{\pm n}$ coefficients}
First, we will derive equation \eqref{lambda2NEW} from Theorem~\ref{main theorem} for the coefficient $\lambda_{+n}^{(2)}$ in the asymptotic expansion \eqref{asymptotic expansion of eigenvalues} of the eigenvalue $n$. We will separate the calculation of this coefficient into several parts, for sake of simplicity and readability, using the explicit formula \eqref{formula za lambda 2}. Let us therefore first evaluate the term  $\langle W_{1/2}^{(2)}v^{(0)},v^{(0)}\rangle$. Using the formula for the eigenvector \eqref{eigenvector} corresponding to the eigenvalue $n$, as well as formula \eqref{A2}  for the differential operator $W_{1/2}^{(2)}$, integrating by parts we get that
\begin{align}\label{part1Lambda2}
\langle W_{1/2}^{(2)} & v^{(0)},v^{(0)}  \rangle
=\int_{0}^{2\pi}[v^{(0)}]^\ast W_{1/2}^{(2)}v^{(0)}dx^1
\\
\nonumber
&=\frac{3n}8 \widehat{(h^2)}{}_{11}(0)
-\frac{n}{8}\tensor{\widehat{k}}{_{11}}(0)
-
\frac{i}{16}\varepsilon_{\beta\gamma1}\sum_{m\in \mathbb{Z}\setminus \{0\}}
m\overline{\widehat{h}_{\alpha\beta}(m)}\widehat{h}_{\alpha\gamma}(m).
\end{align}
We simplified equation \eqref{part1Lambda2} using the Fourier coefficients \eqref{fourier coefficients} and  Parseval's formula
$$\frac{1}{2\pi}\int_{0}^{2\pi}p(x)\overline{q(x)}dx=\sum_{m\in \mathbb{Z}}\widehat{p}(m)\overline{\widehat{q}(m)}.$$
Secondly,  we will evaluate the term $\langle(W_{1/2}^{(1)}-\lambda^{(1)})Q(W_{1/2}^{(1)}-\lambda^{(1)})v^{(0)},v^{(0)}\rangle$.
The pseudoinverse operator \eqref{pseudoinverse operator} corresponding to the operator $W_{1/2}-n I$ is given by
\begin{align}\label{PseudoinverseQ+}
Q_{+n}=\frac{1}{4\pi}
\sum_{m\in\mathbb{Z}\setminus \{n\}}\frac{1}{m-n}
&
\left[
e^{imx^1}
\left(
  \begin{array}{cc}
    1 & 1 \\
    1 & 1 \\
  \end{array}
\right)\int_{0}^{2\pi}e^{-imy^1}(\ \cdot\ )dy^1
\right.
\\
\nonumber
&
\left.+
e^{-imx^1}
\left(
  \begin{array}{cc}
    1 & -1 \\
    -1 & 1 \\
  \end{array}
\right)\int_{0}^{2\pi}e^{imy^1}(\ \cdot\ )dy^1
\right].
\end{align}
Using \eqref{A1}, \eqref{eigenvector} (for eigenvalue $n$) and \eqref{lambda1NEW}, we have that
\begin{align}(W_{1/2}^{(1)}&-\lambda^{(1)}_{+n})v^{(0)}
\nonumber
=
-\frac{n}{4\sqrt{\pi}}
\left(\begin{array}{c}
                       \tensor{h}{_1^1}-i\tensor{h}{_2^1}+\tensor{h}{_3^1}\\
                      \tensor{h}{_1^1}+i\tensor{h}{_2^1}  -\tensor{h}{_3^1}
                    \end{array}
\right)e^{in x^1}\\
\label{A1- lambda1 temp}
&+
\frac{i}{8\sqrt{\pi}}e^{in x^1}\frac{d}{dx^1}
\left(\begin{array}{c}
                      \tensor{h}{_1^1}-i\tensor{h}{_2^1}+\tensor{h}{_3^1}\\
                      \tensor{h}{_1^1}+i\tensor{h}{_2^1}  -\tensor{h}{_3^1}
                    \end{array}
\right)
+\frac{n\tensor{\widehat{h}}{_{11}}(0)}{4\sqrt{\pi}}\left(
      \begin{array}{c}
        1 \\
        1 \\
      \end{array}
    \right)e^{in x^1}.
\end{align}
Using the formula for the pseudoinverse \eqref{PseudoinverseQ+}, we evaluate $Q_{+n}((W_{1/2}^{(1)}-\lambda^{(1)}_{+n})v^{(0)})$  in three parts. First we will act with the pseudoinverse $Q_{+n}$ on the first term on the RHS of equation \eqref{A1- lambda1 temp}.
Using the well known property of the Fourier coefficient that $\widehat{h}(-m)=\overline{\widehat{h}(m)}$, we obtain
\begin{align*}
&Q_{+n}\left(
-\frac{n}{4\sqrt{\pi}}
\left(\begin{array}{c}
\tensor{h}{_1^1}-i\tensor{h}{_2^1}+\tensor{h}{_3^1}\\
\tensor{h}{_1^1}+i\tensor{h}{_2^1}  -\tensor{h}{_3^1}
\end{array}
\right)e^{in x^1}
\right)
=-\frac{n}{8\sqrt{\pi}}
\sum_{m\in\mathbb{Z}\setminus \{n\}}\frac{1}{m-n}
\\&
\left[
e^{imx^1}
\left(
  \begin{array}{cc}
    1 & 1 \\
    1 & 1 \\
  \end{array}
\right)
\frac{1}{2\pi}
\int_{0}^{2\pi}
\left(\begin{array}{c}
                       \tensor{h}{_1^1}-i\tensor{h}{_2^1}+\tensor{h}{_3^1}\\
                      \tensor{h}{_1^1}+i\tensor{h}{_2^1}  -\tensor{h}{_3^1}
                    \end{array}
\right)e^{-i(m-n)y^1}dy^1\right.\\
&+\left.
e^{-imx^1}
\left(
  \begin{array}{cc}
    1 & -1 \\
    -1 & 1 \\
  \end{array}
\right)
\frac{1}{2\pi}
\int_{0}^{2\pi}
\left(\begin{array}{c}
                       \tensor{h}{_1^1}-i\tensor{h}{_2^1}+\tensor{h}{_3^1}\\
                      \tensor{h}{_1^1}+i\tensor{h}{_2^1}  -\tensor{h}{_3^1}
                    \end{array}
\right)e^{-i(-(m+n))y^1}
dy^1
\right]\\
&=-\frac{n}{4\sqrt{\pi}}
\sum_{m\in\mathbb{Z}\setminus \{n\}}\frac{1}{m-n} \\
&
\qquad \left(\begin{array}{c}
                       \widehat{h}_{11}(m-n)e^{imx^1} -i\ \overline{\widehat{h}_{21}(m+n)}e^{-imx^1}
                       +\overline{\widehat{h}_{31}(m+n)}e^{-imx^1}\\
                      \widehat{h}_{11}(m-n)e^{imx^1}+ i\ \overline{\widehat{h}_{21}(m+n)}e^{-imx^1}
                      -\overline{\widehat{h}_{31}(m+n)}e^{-imx^1}
                    \end{array}
\right).
\end{align*}

Secondly, we act with the pseudoinverse $Q_{+n}$ on the second term on the RHS of \eqref{A1- lambda1 temp}.
Integrating by parts, we obtain

\begin{align*}
Q_{+n}&\left(
\frac{i}{8\sqrt{\pi}}e^{in x^1}
\frac{d}{dx^1}
\left(\begin{array}{c}
                       \tensor{h}{_1^1}-i\tensor{h}{_2^1}+\tensor{h}{_3^1}\\
                      \tensor{h}{_1^1}+i\tensor{h}{_2^1}  -\tensor{h}{_3^1}
                    \end{array}
\right)
\right)
=-\frac{1}{16\sqrt{\pi}}
\sum_{m\in\mathbb{Z}\setminus \{n\}}\frac{1}{m-n}
\\
&
\left[
e^{imx^1}
\left(
  \begin{array}{cc}
    1 & 1 \\
    1 & 1 \\
  \end{array}
\right)
\frac{m-n}{2\pi}
\int_{0}^{2\pi}e^{-i(m-n)y^1}
\left(\begin{array}{c}
                       \tensor{h}{_1^1}-i\tensor{h}{_2^1}+\tensor{h}{_3^1}\\
                      \tensor{h}{_1^1}+i\tensor{h}{_2^1}  -\tensor{h}{_3^1}
                    \end{array}
\right)
dy^1\right.\\
&-\left.
e^{-imx^1}
\left(
  \begin{array}{cc}
    1 & -1 \\
    -1 & 1 \\
  \end{array}
\right)
\frac{m+n}{2\pi}
\int_{0}^{2\pi}e^{-i(-(m+n))y^1}
\left(\begin{array}{c}
                       \tensor{h}{_1^1}-i\tensor{h}{_2^1}+\tensor{h}{_3^1}\\
                      \tensor{h}{_1^1}+i\tensor{h}{_2^1}  -\tensor{h}{_3^1}
                    \end{array}
\right)dy^1
\right]\\
&=-\frac{1}{8\sqrt{\pi}}
\sum_{m\in\mathbb{Z}\setminus \{n\}}
\left(\begin{array}{c}
                       \widehat{h}_{11}(m-n)
                       \\
                      \widehat{h}_{11}(m-n)
                    \end{array}
\right)e^{imx^1}\\
&-\frac{1}{8\sqrt{\pi}}
\sum_{m\in\mathbb{Z}\setminus \{n\}}\frac{m+n}{m-n}\left(\begin{array}{c}
                       i\overline{\widehat{h}_{21}(m+n)}
                       -\overline{\widehat{h}_{31}(m+n)}\\
                      -i\overline{\widehat{h}_{21}(m+n)}
                      +\overline{\widehat{h}_{31}(m+n)}
                    \end{array}
\right)e^{-imx^1}
.
\end{align*}
Thirdly, we act with $Q_{+n}$ on the third and final term on the RHS of \eqref{A1- lambda1 temp}
\begin{align}\label{suma1}
Q_{+n}&\left(\frac{\tensor{\widehat{h}}{_{11}}(0)}{4\sqrt{\pi}}\left(
      \begin{array}{c}
        1 \\
        1 \\
      \end{array}
    \right)e^{in x^1}\right)= \\
\nonumber
&    \frac{\tensor{\widehat{h}}{_{11}}(0)}{16\pi\sqrt{\pi}}
\sum_{m\in\mathbb{Z}\setminus \{n\}}\frac{1}{m-n}
\left[
e^{imx^1}
\left(
  \begin{array}{cc}
    1 & 1 \\
    1 & 1 \\
  \end{array}
\right)\int_{0}^{2\pi}\left(
      \begin{array}{c}
        1 \\
        1 \\
      \end{array}
    \right)e^{i(n-m)y^1}dy^1\right.\\
\nonumber
&
+\left.
e^{-imx^1}
\left(
  \begin{array}{cc}
    1 & -1 \\
    -1 & 1 \\
  \end{array}
\right)\int_{0}^{2\pi}\left(
      \begin{array}{c}
        1 \\
        1 \\
      \end{array} \right)
    e^{i(m+n)y^1}dy^1
\right].
\end{align}
For $m\in \mathbb{Z}$ we have that
$$\int_{0}^{2\pi}e^{i(n-m)y^1}dy^1=
\left\{
\begin{array}{cc}
  2\pi, & m=n ,\\
  0, & m\neq n,
\end{array}
\right.
$$
$$
\int_{0}^{2\pi}e^{i(m+n)y^1}dy^1=
\left\{
\begin{array}{cc}
  2\pi, & m=-n,\\
  0, & m\neq -n.
\end{array}
\right.
$$
We get that  the sum  \eqref{suma1} only makes sense if $m=-n$. Hence
\begin{align*}
&Q_{+n}\left(\frac{\tensor{\widehat{h}}{_{11}}(0)}{4\sqrt{\pi}}\left(
      \begin{array}{c}
        1 \\
        1 \\
      \end{array}
    \right)e^{in x^1}\right)=
-\frac{\tensor{\widehat{h}}{_{11}}(0)}{32\pi\sqrt{\pi}}
e^{in x^1}
\left(
  \begin{array}{cc}
    1 & -1 \\
    -1 & 1 \\
  \end{array}
\right)\left(
      \begin{array}{c}
        2\pi \\
        2\pi \\
      \end{array} \right)
     =0.
\end{align*}
Putting the above  calculations together, we get that
\begin{align}\label{tempQAlambda}
Q_{+n}&((A^{(1)}-\lambda^{(1)}_{+n})v^{(0)})=
-\frac{1}{8\sqrt{\pi}}
\sum_{m\in\mathbb{Z}\setminus \{n\}}
\frac{m+n}{m-n}
\left(\begin{array}{c}
                       \widehat{h}_{11}(m-n)\\
                      \widehat{h}_{11}(m-n)
                    \end{array}
\right)e^{imx^1}
\\ \nonumber &
-\frac{1}{8\sqrt{\pi}}
\sum_{m\in\mathbb{Z}\setminus \{n\}}
\left(\begin{array}{c}
                       i\overline{\widehat{h}_{21}(m+n)}
                       -\overline{\widehat{h}_{31}(m+n)}\\
                      -i\overline{\widehat{h}_{21}(m+n)}
                      +\overline{\widehat{h}_{31}(m+n)}
                    \end{array}
\right)e^{-imx^1}.
\end{align}
Calculating $\frac{d(Q_{+}((A^{(1)}-\lambda^{(1)}_{+})v^{(0)}))}{dx^1}$ and
using equations \eqref{A1}, \eqref{lambda1NEW} and \eqref{tempQAlambda}, we get the term $(W_{1/2}^{(1)}-\lambda_{+}^{(1)})Q_{+}((W_{1/2}^{(1)}-\lambda^{(1)}_{+})v^{(0)})$
and finally, the second term
$\langle(W_{1/2}^{(1)}-\lambda_{+}^{(1)})Q_{+}(W_{1/2}^{(1)}-\lambda^{(1)}_{+})v^{(0)},v^{(0)}\rangle$,
using \eqref{scalar product}, becomes

\begin{align}\label{part2Lambda2}
\langle(W_{1/2}^{(1)}&
-\lambda_{+n}^{(1)})Q_{+n}((W_{1/2}^{(1)}
-\lambda^{(1)}_{+n})v^{(0)}),v^{(0)}\rangle\\
\nonumber
=&
\frac{1}{16}\sum_{m\in \mathbb{Z}\setminus \{n\}} \frac{1}{m-n}(m+n)^2\ \widehat{h}_{11}(m-n)\overline{\widehat{h}_{11}(m-n)}\\
\nonumber
+&\frac{1}{16}\sum_{m\in \mathbb{Z}\setminus \{n\}}(m-n)\
\widehat{h}_{31}(m+n)
\left(\overline{\widehat{h}_{31}(m+n)}
-i\overline{\widehat{h}_{21}(m+n)}\right)\\
\nonumber
+&\frac{1}{16}\sum_{m\in \mathbb{Z}\setminus \{n\}}i(m-n)\
\widehat{h}_{21}(m+n)
\left(\overline{\widehat{h}_{31}(m+n)}
-i\overline{\widehat{h}_{21}(m+n)}\right).
\end{align}
Combining  equations \eqref{part1Lambda2} and \eqref{part2Lambda2}, we get the formula \eqref{formula za lambda 2}, for the coefficient $\lambda_{+}^{(2)}$ is given by \eqref{lambda2NEW}.
\begin{remark}
Using the eigenvector \eqref{eigenvector} corresponding to the eigenvalue $-n$ of the massless Dirac operator and the pseudoinverse operator \eqref{pseudoinverse operator} (corresponding to $-n$) of the operator $W_{1/2}+n I$, analogously to the above calculations performed for the eigenvector  \eqref{eigenvector} corresponding to the eigenvalue $n$, we get the formula \eqref{formula za lambda 2}, for the coefficient $\lambda_{-}^{(2)}$.
\end{remark}

\section{Numerical analysis of our results}
This appendix deals with the operator \eqref{massless dirac} and here we numerically analyze its spectrum. Using Galerkin's method (see e.g. \cite{grossmann2007numerical}), we discretize the eigenvalue problem of this operator.

Consider the  $2m+1$  eigenvalues
$\lambda_i=i,\  (i=0,\pm 1,\ldots,\pm m)$
of the unperturbed massless Dirac operator on half-densities $W_{1/2}(0)$. Each eigenvalue $\lambda_i$ 
  has  multiplicity two and the corresponding eigenvectors $v_i(x^1)$ and $w_i(x^1)$,  $(i=0, \pm 1, \ldots, \pm m)$  are given  by
\eqref{eigenvector} and \eqref{eigenvector Cv}. Hence, we have that
\begin{align}
\label{eigenvalueProblemOfUnperturbedOperator1}
W_{1/2}(0)v_i(x^1)&=\lambda_i v_i(x^1),\\
\label{eigenvalueProblemOfUnperturbedOperator2}
W_{1/2}(0)w_i(x^1)&=\lambda_i w_i(x^1),
\end{align}
where $i=0,\pm 1,\ldots,\pm m. $
The eigenvectors $v_i(x^1)$ and $w_i(x^1)$ are orthonormal with respect to the inner product \eqref{scalar product}, i.e.
\begin{equation}\label{eigenvectors are orthonormal}
\langle v_i, v_j\rangle=\langle w_i, w_j\rangle=
\delta_{ij}, \ \ \
\langle v_i,w_j\rangle=\langle w_i,v_j\rangle=0, \ \ \  (i,j=0,\pm 1,\ldots,\pm m).
\end{equation}
According to \eqref{eigenvectors are orthonormal}, from equations \eqref{eigenvalueProblemOfUnperturbedOperator1} and \eqref{eigenvalueProblemOfUnperturbedOperator2}, for $i,j=0,\pm 1,\ldots,\pm m$ we have that
\begin{align*}
\lambda_i&=\langle W_{1/2}(0)v_i(x^1),v_i(x^1)\rangle
=\langle W_{1/2}(0)w_i(x^1),w_i(x^1)\rangle
\end{align*}
and
\begin{align*}
\langle W_{1/2}(0)v_i(x^1),w_j(x^1)\rangle &=\langle W_{1/2}(0)v_j(x^1),w_i(x^1)\rangle=0, \\
\langle W_{1/2}(0)w_i(x^1),v_j(x^1)\rangle &=\langle W_{1/2}(0)w_j(x^1),v_i(x^1)\rangle =0.
\end{align*}
 Let us now construct the matrices
\begin{equation}\label{blok matrica Hmn}
H_{i,j}:=\left(
           \begin{array}{cc}
             \langle W_{1/2}(0)v_{i},v_{j}\rangle & \langle W_{1/2}(0)v_{i},w_{j}\rangle \\
             \langle W_{1/2}(0)w_{i},v_{j}\rangle & \langle W_{1/2}(0)w_{i},w_{j} \rangle\\
           \end{array}
         \right)
\end{equation}
where $i,j=0,\pm 1,\ldots,\pm m.$ Using matrices (\ref{blok matrica Hmn}), we can construct the block matrix $H$ as follows
\begin{equation*}
H:=
\left(
  \begin{array}{ccccccc}
    H_{-m,m} &  &  & H_{0,m} &  &  & H_{m,m} \\
    &\ddots&&\vdots&&\iddots & \\
     &  &  & H_{0,1} &  &  &  \\
     & \cdots & H_{-1,0} & H_{0,0} & H_{1,0} & \cdots &  \\
     &  &  & H_{0,-1} &  &  &  \\
    &\iddots &&\vdots&&\ddots& \\
    H_{-m,-m} & &  & H_{0,-m} &  &  & H_{m,-m} \\
  \end{array}
\right).
\end{equation*}
The  matrix $H$ is a quadratic matrix of order $2(2m+1)$ and by construction it is a Hermitian matrix.
The eigenvalues of the matrix $H$ are $\lambda=0,\pm 1,\ldots ,\pm m$ and each eigenvalue has multiplicity two. Therefore, the analysis of the eigenvalue problem \eqref{eigenvalueProblemOfUnperturbedOperator1}-\eqref{eigenvalueProblemOfUnperturbedOperator2} and the analysis of the eigenvalue problem of the matrix $H$ are equivalent.

Now we choose to consider the matrix $H(\upepsilon)$ with the perturbed massless Dirac operator $W_{1/2}(\upepsilon)$ instead of the unperturbed operator $W_{1/2}(0)$. Then the matrix $H(\upepsilon)$ is
a Hermitian matrix whose entries depend on the parameter $\upepsilon$ and $H(0)=H$. Using the perturbation process described in \cite{mccartin2003pseudoinverse} for  perturbed Hermitian  matrices, we can get the asymptotic expansions of the eigenvalues of the perturbed matrix $H(\upepsilon)$ and specially the asymptotic expansions of the eigenvalues $\lambda=\pm 1$. The eigenvalues of the matrix $H(\upepsilon)$ will converge to the eigenvalues of the matrix $H$ as $\upepsilon\rightarrow 0.$
The examination of the spectrum of the perturbed massless Dirac operator will reduce to the examination of the spectrum of the Hermitian matrix $H(\upepsilon)$. The numerical calculations were performed using \emph{Wolfram Mathematica}.

\begin{example}
Consider the coframe
\begin{equation}\label{tempCoframe}
\tensor{e}{^j_\alpha}=\tensor{\delta}{^j_\alpha}+
\upepsilon
\left(
  \begin{array}{ccc}
    0 & 0 & 0 \\
    0 & \cos x^1 & \sin x^1 \\
    0 & \sin x^1 & -\cos x^1 \\
  \end{array}
\right).
\end{equation} which was also considered in \cite{downes2013spectral}.
The explicit formula for the perturbed massless Dirac operator corresponding to the coframe \eqref{tempCoframe}  reads
\begin{equation}\label{temp4}
W(\upepsilon)=-i
\left(
  \begin{array}{cc}
    0 & 1 \\
    1 & 0 \\
  \end{array}
\right)\frac{d}{dx^1}-
\frac{\upepsilon^2}{2(1-\upepsilon^2)}I.
\end{equation}
The eigenvalues of the operator \eqref{temp4} are explicitly given by
\begin{equation}\label{tempLambda}
\lambda_n(\upepsilon)=n-\frac{\upepsilon^2}{2(1-\upepsilon^2)}=n-\frac{\upepsilon^2}{2}-\frac{\upepsilon^4}{2}+O(\upepsilon^6),\ \ n\in \mathbb{Z}
\end{equation}
and all eigenvalues have multiplicity two.

Now we will use the coframe \eqref{tempCoframe} to analyze the spectrum of the massless Dirac operator using the Galerkin method described above in order to numerically confirm these results. We explicitly constructed the matrix $H(\upepsilon)$ of order $102\times 102$
and numerically analyzed the part of its spectrum.
The  eigenvalues $0,\pm 1,\pm 2$  of the matrix $H(\upepsilon)$ are perturbed as follows
\begin{center}
\begin{tabular}{|l|c|c|c|c|c|}
  \hline
   & $-2$ & $-1$ & $0$ & $1$ & $2$ \\ \hline\hline
   $\upepsilon=0.2$ & -2.02083&-1.02083&-0.0208333&0.979167&1.97917 \\ \hline
  $ \upepsilon=0.1$ & -2.00505&-1.00505&-0.00505051& 0.994949&1.994950 \\ \hline
   $\upepsilon=0.01$ &-2.00005&-1.00005&-0.000050005& 0.99995& 1.99995 \\
  \hline
\end{tabular}
\end{center}
and each eigenvalue has  multiplicity two.
Analyzing the data given in the above table we see that for this choice of the coframe the spectral symmetry of the matrix $H(\upepsilon)$ is broken and consequently we obtain spectral asymmetry of the massless Dirac operator in the axisymmetric case.

Using the perturbation process for the matrices with double eigenvalues described in \cite{mccartin2003pseudoinverse}, we get that the asymptotic formulae for the eigenvalues $\pm 1$ are given by
\begin{align*}
\lambda_{+}(\upepsilon)&=1-\frac12\upepsilon^2+O(\upepsilon^3),\\
\lambda_{-}(\upepsilon)&=-1-\frac12\upepsilon^2+O(\upepsilon^3).
\end{align*}
which is in accordance with \eqref{tempLambda}.
\end{example}
Now, we will consider the coframe which is not symmetric  to show that in this case  it is also possible to obtain spectral asymmetry.
\begin{example}
Consider the coframe
\begin{equation*}
\tensor{e}{^j_\alpha}=\tensor{\delta}{^j_\alpha}+
\upepsilon
\left(\!
  \begin{array}{ccc}
    0 & \cos x^1-\cos 2x^1+\cos 3x^1 & \sin x^1+\sin 2x^1-\sin 3x^1 \\
    0 & 0 & 0 \\
    0 & 0 & 0 \\
  \end{array}\!
\right)\!.
\end{equation*}
Analyzing the spectrum of the matrix  $H(\upepsilon)$ of order $102\times 102$  we get that the  eigenvalues $0,\pm 1,\pm 2$  of the matrix  are perturbed as follows
\begin{center}
\begin{tabular}{|l|c|c|c|c|c|}
  \hline
   & $-2$ & $-1$ & $0$ & $1$ & $2$ \\ \hline\hline
   $\upepsilon=0.2$ & -2.10913& -1.05372& 0.00169489& 1.0571& 2.11252 \\ \hline
  $\upepsilon=0.1$ & -2.02923& -1.01456& 0.000119453& 1.0148& 2.02947 \\ \hline
  $\upepsilon=0.01$ &-2.0003&-1.00015&$1.24941\times 10^{-8}$& 1.00015& 2.0003 \\
  \hline
\end{tabular}
\end{center}
and each eigenvalue has  multiplicity two. Analyzing the data given in the above table  we see that the  spectral symmetry is broken.
Using the method described in \cite{mccartin2003pseudoinverse}, we obtain that the  asymptotic formulae for the eigenvalues $\pm 1$ are given by
\begin{align*}
\lambda_{+}(\upepsilon)&=1+\frac32\upepsilon^2-\frac{17}{8}\upepsilon^4+O(\upepsilon^5),\\
\lambda_{-}(\upepsilon)&=-1-\frac32\upepsilon^2+\frac{37}{8}\upepsilon^4+O(\upepsilon^5),
\end{align*}
hence we see that spectral asymmetry is achieved in the quartic term.
\end{example}

\bibliographystyle{amsplain}
\bibliography{Analyzing-the-spectral-asymmetry-of-the-massless-Dirac-operator-on-the-3-torus-v2}
\label{lastpage}
\end{document}